\documentclass{amsart}
\usepackage[margin=1in]{geometry}
\usepackage{mathrsfs,nicefrac,graphicx,float,cancel,tikz,esdiff,mathtools,bm}
\usepackage{lipsum}
\usepackage{float,graphicx,subcaption}
\usepackage{enumitem,hyperref,nameref}
\usepackage{subfiles}
\usepackage[toc,page]{appendix}
\usepackage[utf8]{inputenc}

\newtheorem{theorem}{Theorem}
\newtheorem{lemma}[theorem]{Lemma}
\newtheorem{proposition}{Proposition}
\newtheorem{definition}{Definition}

\newtheorem{problem}[theorem]{Problem}
\def\thmheadbrackets#1#2#3{%
	\thmname{#1}\thmnumber{\@ifnotempty{#1}{ }\@upn{#2}}%
	\thmnote{ {\the\thm@notefont[#3]}}}
\makeatother

\newtheoremstyle{brakets}
{}
{}
{\itshape}
{}
{\bfseries}
{.}
{ }
{\thmheadbrackets{#1}{#2}{#3}}

\theoremstyle{brakets}

\theoremstyle{remark}
\newtheorem*{remark}{Remark}

\newcommand{\RR}{\mathbb R}
\newcommand{\ZZ}{\mathbb Z} 
\newcommand{\norm}[1]{\ensuremath{\left\lVert #1 \right\rVert}}
\newcommand{\abs}[1]{\ensuremath{\left\lvert #1 \right\rvert}}
\newcommand{\innerpdt}[1]{\ensuremath{\left<#1\right>}}

\let\phi\varphi
\let\varphi\phi

\let\bar\overline
\let\overline\bar

\newcommand{\R}{\mathbb{R}}

\newcommand{\C}{\mathbb{C}}

\newcommand{\be}{\begin{enumerate}}
\newcommand{\ee}{\end{enumerate}}
\newcommand{\bi}{\begin{itemize}}
\newcommand{\ei}{\end{itemize}}

\newcommand{\ba}{\begin{array}}
\newcommand{\ea}{\end{array}}
\newcommand{\bmat}{\left[\begin{array}}
\newcommand{\emat}{\end{array}\right]}
\newcommand{\bt}{\begin{thm}}
\newcommand{\et}{\end{thm}}
\newcommand{\bp}{\begin{proof}}
\newcommand{\ep}{\end{proof}}
\newcommand{\bprop}{\begin{prop}}
\newcommand{\eprop}{\end{prop}}
\newcommand{\bl}{\begin{lemma}}
\newcommand{\el}{\end{lemma}}
\newcommand{\bc}{\begin{cor}}
\newcommand{\ec}{\end{cor}}
\newcommand{\bd}{\begin{defn}}
\newcommand{\ed}{\end{defn}}

\newcommand{\vep}{\varepsilon}

\begin{document}
\title{Stationary Curves Under the M\"obius-Plateau Energy}
\author{Max Lipton}
\address{Department of Mathematics, Massachusetts Institute of Technology, MA 02139}
\email{liptonm@mit.edu}

\author{Gokul G. Nair}
\address{Department of Mathematics, Rutgers University, NJ 08854}
\email{gokul.nair@rutgers.edu}

\begin{abstract}
Plateau problems with elastic boundary energies have been of recent theoretical and applied interest. However, strong assumptions have to be made to avoid self-intersections of the boundary curve during energy minimization. We introduce a class of Plateau problems for boundaries with self-repulsive energies that obviates self-contact in energy minimization problems. For the self-repulsive energy, we choose the M\"obius Energy introduced by O'Hara due to its myriad regularity properties shown by Freedman et al. We first prove an existence theorem for this \textit{M\"obius-Plateau} problem in the class of closed Lipschitz curves of a given irreducible knot-type spanned by immersed discs. We then turn our attention to M\"obius-Plateau variations of helicoidal strips, which are classified as ``screw-like" or ``ribbon-like" based on the signs of the radii of the boundary helices. By analyzing the Euler-Lagrange equations, we show that screw-like solutions are plentiful, whilst ribbon-like solutions impose strong constraints on their parameters: they must have high frequency (equivalently, low pitch), thin width in comparison to the frequency, and remain close to the axis.
\end{abstract}

\maketitle

\section{Introduction}
The problem of finding a surface of minimal area bounded by a fixed closed curve, also known as \textit{Plateau's problem}, has captivated mathematicians for over two centuries. A relevant extension to this problem is when the boundary possesses its own energy and is allowed to deform. Since most experimental realizations of Plateau's problem involve films formed by dipping closed loops of wire in soapy water, it is natural to consider boundary curves with elastic energy, defined as the integral of squared curvature. This problem, known now as the \textit{Euler-Plateau} problem has received considerable attention in recent years~\cite{bernatzki2001minimal,giomi2012minimal,gruber21,chen2014stability}. In~\cite{bernatzki2001minimal}, the authors prove the existence of loop configurations that minimize a combination of the elastic energy of curves with preferred curvature and the area of a minimal surface spanning the curve. However, they require strong assumptions to avoid self-contact of the boundary curve during energy minimization. The paper \cite{giusteri2017solution} models the boundary curve as a Kirchhoff rod and find minimizers whose boundary curves can touch without interpenetration of matter.

In view of this, we propose an interesting variation on the above theme--- minimal surfaces bounded by self-repulsive curves. Competition between the self-repulsive force trying to enlarge the bounding curve and surface tension of the spanning surface trying to shrink it leads to interesting physics in this problem. Knots charged with regularized Coulomb-like self-interaction energies were considered originally by 
O'Hara~\cite{o1991energy,o1992family,o2003energy} and studied extensively in~\cite{freedman1994mobius}. In particular, it is shown in~\cite{freedman1994mobius} that this energy is M\"obius invariant, and is hence referred to as the \textit{M\"obius} energy. Furthermore, \cite{freedman1994mobius} establishes the existence of loops minimizing the M\"obius energy in fixed irreducible knot classes. In this work, we consider an energy that combines the M\"obius self-repulsive energy of a closed loop and the area of a minimal disk spanning that loop and prove the existence of minimizing configurations for this energy. We apply well-known results on M\"obius minimizers stated in \cite{freedman1994mobius}, but lower semicontnuity of the Plateau energy, while assumed in the Euler-Plateau literature (e.g. \cite{bernatzki2001minimal}), does not have a written proof to the authors' knowledge. One might be misled into believing that the sequence of Plateau solutions associated to a given uniformly convergent sequence of boundary curves also converges, but even when the solutions catastrophically jump when passing to the limit, the area is still lower semicontinuous, masking the catastrophe.

And so, even though minimizers of the M\"obius-Plateau energy exist, the competition between the two constituent energies can complicate the critical points themselves. We then turn our attention to the corresponding Euler-Lagrange equations and explore an example of this phenomenon, even though our prior existence result does not apply. In particular, we consider helicoidal strips--- regions bounded between two helices on a common helicoid as candidate solutions. These configurations can be classified as helicoidal screws or ribbons depending on the signs of the radii of the bounding helices. By working directly with the oscillatory variational integrals to derive nontrivial necessary conditions, we will see that while screw-like critical points are abundant, ribbon-like solutions are strongly constrained: they must have high-frequency (low pitch), thin width and remain close to the axis. Many of the calculations in this section were performed with the assistance of a computer algebra system for the sake of practicality.

The outline of this paper is as follows: In Section~\ref{sec:problem-formulation} we recall the definition and properties of the M\"obius energy and formulate the M\"obius-Plateau energy minimization problem. In Section~\ref{sec:energy-minimization}, we show compactness properties and lower semicontinuity of the energy and establish the existence of disk-like energy minimizing configurations for simple closed curves in a given irreducible knot type. In Section~\ref{sec:stationary-helix-pairs} we turn our attention to helicoidal solutions of the Euler-Lagrange equations.  

\section{Problem Formulation}\label{sec:problem-formulation}
Denote the closed unit disk by $D\subset\RR^2$ and let $\gamma:\partial D\rightarrow\RR^3$ be a rectifiable simple closed curve. We define the \textit{M\"obius energy}~\cite{o2003energy,freedman1994mobius} as
\begin{align}\label{def:mobius-energy}
    E_M[\gamma]=\iint_{\partial D\times\partial D}\left[\frac{1}{\abs{\gamma(y)-\gamma(x)}^2}-\frac{1}{D(\gamma(y),\gamma(x))^2}\right]\abs{\dot{\gamma}(x)}\abs{\dot{\gamma}(y)}d yd x,
\end{align}
where $D(\gamma(x),\gamma(y))$ denotes the geodesic distance between $\gamma(x)$ and $\gamma(y)$ along the curve $\gamma$. The M\"obius energy blows up in self intersections, so the gradient descent will not encounter boundary catastrophes, an issue which can occur in Euler-Plateau problems requiring strong assumptions to avoid \cite{bernatzki2001minimal}. Furthermore, minimizers of the M\"obius energy enjoy nice regularity properties. Freedman, He, and Wang proved the minimizers are $C^{1,1}$ \cite{freedman1994mobius}, a result later strengthened to $C^\infty$ regularity by He \cite{he2000}, and finally to analyticity by Blatt and Vorderobermeier \cite{blatt18}.

For a fixed curve $\gamma$, we define the following class of parametrized discs that span the curve,
\begin{align*}
    \mathcal{D}_\gamma:=\{&u:D\rightarrow\RR^3|~u\in W^{1,2}(D)\cap C^0(\bar{D}),~u:\partial D\rightarrow\gamma\text{ is monotone and onto}\}
\end{align*}
The area of the image of $u\in\mathcal{D}_\gamma$ can be computed as
\begin{align}\label{def:area}
    \mathrm{Area}(u)=\int_D \left(\abs{u_x}^2\abs{u_y}^2-\innerpdt{u_x,u_y}^2\right)^{\frac{1}{2}}d xd y
\end{align}
where $u_x$ and $u_y$ are the derivatives of $u$ in $x$ and $y$. We then have the following definition of the \textit{Plateau energy},
\begin{align}\label{eqn:plateau-energy}
    E_P[\gamma]=\inf_{u\in\mathcal{D}_\gamma}\mathrm{Area}(u).
\end{align}
Note the Plateau energy is defined on curves, not surfaces. Furthermore, we restrict the topological class of the candidate solutions to a disk. Our problem of interest can be stated as follows,
\begin{problem}
    For fixed constants $\alpha, \beta \geq 0$, find a closed curve $\gamma$ of prescribed length and knot type which minimizes the \textbf{M\"obius-Plateau energy},
    \begin{align}\label{def:mobius-plateau-energy}
        E[\gamma]= \alpha E_M[\gamma]+ \beta E_P[\gamma].
    \end{align}
\end{problem}

\section{Energy Minimization}\label{sec:energy-minimization}
In this section, we fix $\alpha = \beta = 1$. Recall that a knot type refers to an equivalence class of rectifiable simple closed curves in $\R^3$ with respect to ambient isotopy. A knot type is irreducible if a representative cannot be realized as the connect sum of two nontrivially knotted curves. By convention, we will also say the unknot is irreducible. We restrict to irreducible knot types in order to apply the results of \cite{freedman1994mobius}. The existence of M\"obius minimizers for composite knot types remains open due to the possibility of ``pull-tight" singularities in the gradient descent. 

We define the following admissible class of closed curves for some irreducible knot type $K$,
\begin{align*}
    \mathcal{A}=\{&\gamma\in C^{0,1}(\partial D,\RR^3)|~\gamma(\partial D)\text{ is simple closed},\\
    &\gamma\text{ is of knot type } K,~\mathrm{Length}(\gamma) = 2\pi,~\gamma(0)=0\} 
\end{align*}
Here, $C^{0,1}(\partial D,\RR^3)$ denotes Lipschitz maps from $\partial D\cong S^1$ to $\RR^3$. The constraints $\gamma(0)=0$ (written with an abuse of notation) and $\mathrm{Length}(\gamma) = 2\pi$ together ensure that the images of the curves in $\mathcal{A}$ are contained in some ball of fixed radius. This compactness condition allows for lower semicontinuity with respect to uniform convergence to imply the existence of a minimizer.

We also state some key results that are useful in the proof of the existence theorem.
\begin{lemma}[\cite{freedman1994mobius} Lemma 1.2]\label{lemma:freedman-lemma-1.2}
    Let $\gamma:X\rightarrow\RR^3$ be a rectifiable curve in $\RR^3$ parametrized by arc length. If $E_M[\gamma]$ is finite, then $\gamma$ is $C$ bi-Lipschitz with constant $C$ depending only on $E_M[\gamma]$. Furthermore, $C \to 1$ as $E_M[\gamma] \to 0$. 
\end{lemma}

As noted in \cite{freedman1994mobius}, the M\"obius energy enjoys a property akin to absolute continuity: subarcs of a curve with finite M\"obius energy can be chosen to have arbitrarily small energy provided that their length is sufficiently small. Since the M\"obius energy is small, the (inverse) intrinsic and extrinsic distances are close in $L^2$, and so the subarcs are close to being isometrically embedded in $\R^3$.

\begin{lemma}[\cite{freedman1994mobius} Corollary 1.3]\label{lemma:freedman-cor-1.3} With $\gamma$ as above such that $E_M[\gamma]$ is finite, we have that for any $\varepsilon > 0$, there exists $\delta = \delta(\varepsilon, \gamma) > 0$ such that $\gamma$ is a $(1+\varepsilon)$ bi-Lipschitz embedding on subarcs of $\gamma$ with length less than $\delta$.
\end{lemma}

The $(1 + \varepsilon)$ bi-Lipschitz property in the prior lemma is a weaker notion of an isometric embedding. In proving the existence of minimizers for the M\"obius-Plateau energy, we will need to ensure that arc length is preserved in the limit of a minimizing sequence. Were we afforded the luxury of our curves being isometric embeddings of fixed length, length preservation in the limit would trivially follow. In general, total length is not preserved in sequences of curves, even in a uniformly converging limit. In particular, the length could shorten. However, we will soon see this weaker property is sufficient for the preservation of total length.

\begin{lemma}[\cite{freedman1994mobius} Lemma 4.2]\label{lemma:freedman-lemma-4.2}
    Let $\gamma_i:\RR/L_i\ZZ\rightarrow\RR^3$ be a sequence of rectifiable simple closed curves of uniformly bounded energy, $E_M$. Assume that the curves are all parametrized by arc-length and that $L=\lim_{i\rightarrow\infty}L_i>0$ exists. If $\gamma_i(0)$ is a bounded sequence of points, then there is a subsequence $\gamma_{i_k}$ of $\gamma_i$ which converges uniformly to a rectifiable simple curve $\gamma:\RR/L\ZZ\rightarrow\RR^3$. Moreover, $E_M$ is sequentially lower semicontinuous under this convergence, i.e.,
    \begin{align*}
        E_M[\gamma]\leq\liminf_{i_k\rightarrow\infty}E_M[\gamma_{i_k}]
    \end{align*}
\end{lemma}

\begin{theorem}[\cite{freedman1994mobius} Theorem 4.3]\label{thm:freedman-thm-4.3}
    Let $K$ be an irreducible knot. There exists a simple loop $\gamma_K: \RR/L\ZZ\rightarrow\RR^3$ with knot-type $K$ such that $E_M[\gamma_K]\leq E_M[\gamma]$ for any other simple closed loop $\gamma:\RR/L\ZZ\rightarrow\RR^3$ of the same knot-type.
\end{theorem}

We begin the proof of existence by considering a minimizing sequence of $E$, $\{\gamma_k\}\subset\mathcal A$ i.e.,
\begin{align*}
    \lim_{k\rightarrow\infty}E[\gamma_k]=\inf_{\gamma\in\mathcal{A}}E[\gamma].
\end{align*}
This gives us a uniform bound on $E_M[\gamma_k]$, therefore, from Lemma~\ref{lemma:freedman-lemma-4.2} we have a subsequence (for which we will use the same labelling) that converges uniformly to a simple closed $C^{0,1}$ curve $\gamma_\infty: \partial D \to \RR^3$. Furthermore, by applying the remark after Theorem~\ref{thm:freedman-thm-4.3} we can conclude that $\gamma_\infty$ has the same knot-type $K$, even though the M\"obius energy of the limit is not necessarily the infimum. It is also straightforward to conclude that $\gamma_\infty(0)=0$ and apply Fatou's Lemma to see
\begin{align*}
    \mathrm{Length}(\gamma_\infty)&=\int_{0}^{2\pi}\abs{\dot{\gamma}_\infty(s)}d s\leq\liminf_{k\rightarrow\infty}\int_{0}^{2\pi}\abs{\dot{\gamma}_k(s)}d s=\liminf_{k\rightarrow\infty}~\mathrm{Length}(\gamma_k) =  2\pi.
\end{align*}

We must still show $\text{Length}(\gamma_\infty) \geq 2\pi$. Let $\varepsilon > 0$. Applying Lemma~\ref{lemma:freedman-cor-1.3} we have that each $\gamma_i$ is a $(1+\varepsilon)$ bi-Lipschitz embedding on subarcs of length $\leq \delta = \delta(i,\varepsilon)$. As the sequence converges locally uniformly, $\gamma_\infty$ is also a $(1+\varepsilon)$ bi-Lipschitz embedding, and we can assume this $\delta$ is independent of $i$ and also works for $\gamma_\infty$. Now take a partition $0 = s_0 < s_1 < \dots < s_N = 2\pi$ with $|s_{i+1} - s_i| < \delta$ and $N$ large enough so that $N\delta > 2\pi$. By the bi-Lipschitz property, we have $|\gamma_\infty(s_{i+1}) - \gamma_\infty(s_i)| \geq (1+\varepsilon)^{-1} |s_{i+1} - s_i| > (1 + \varepsilon)^{-1}\delta$. Summing this inequality over $i$ yields $\text{Length}(\gamma_\infty) \geq (1 + \varepsilon)^{-1}N\delta > (1 + \varepsilon)^{-1}2\pi$. As $\varepsilon$ was arbitrary, we conclude $\text{Length}(\gamma_\infty) \geq 2\pi$, as desired, so therefore $\gamma_\infty \in \mathcal{A}$.

We would like to now show that $E$ is weakly lower semicontinuous with respect to uniform convergence $\gamma_k\rightarrow\gamma_\infty$. We already have lower semicontinuity in $E_M$. So it remains to show that $E_P$ is lower semicontinuous, i.e.
\begin{align}\label{eqn:E_P-lower-semicontinuity1}
    E_P[\gamma_\infty]\leq\liminf_{k \to \infty} E_P[\gamma_k]
\end{align}

\subsection{Lower Semicontinuity of \texorpdfstring{$E_P$}{Plateau Energy}}
We define the Dirichlet energy of a mapping as 
\begin{align*}\label{def:dirichlet-energy}
    \mathrm{e}(u) = \int_D \abs{\nabla u}^2d xd y.
\end{align*}
\begin{lemma}\label{lemma:area-dirichlet-equivalence}
    Let $\gamma\in\mathcal{A}$, then
    \begin{align*}
        \inf_{u\in\mathcal{D}_\gamma}\mathrm{Area}(u)=\inf_{u\in\mathcal{D}_\gamma}\mathrm{e}(u)
    \end{align*}
\end{lemma}
\begin{proof}
    Refer to~\cite{colding2011course} Lemma 4.4.
\end{proof}

First we state the following version of the classical Plateau's problem. See ~\cite{douglas1931solution,colding2011course} for the Lipschitz and piecewise $C^1$ cases, and \cite{simongmt} for the rectifiable case.
\begin{theorem}\label{theorem:colding-theorem-4.1}
    Given a rectifiable simple closed curve $\Gamma\subset\RR^3$, there exists a map $u:D\rightarrow\RR^3$ such that
    \begin{enumerate}
        \item $u:\partial D\rightarrow\Gamma$ is monotone and onto.
        \item $u\in C^0(\bar D)\cap W^{1,2}(D)$
        \item The image of $u$ minimizes area among all maps from disks with boundary $\Gamma$.
    \end{enumerate}
\end{theorem}
Applying Theorem~\ref{theorem:colding-theorem-4.1} to each of the simple closed curves $\gamma_k(\partial D)$ $k=1,2,...,\infty$ and using Lemma~\ref{lemma:area-dirichlet-equivalence}, we get a sequence of Plateau solutions $u_k\in\mathcal{D}_{\gamma_k}$. That is, for all $k$,
\begin{align*}
    \mathrm{e}(u_k)=\mathrm{Area}(u_k)=\inf_{u\in\mathcal{D}_{\gamma_k}}\mathrm{Area}(u)=E_P[\gamma_k].
\end{align*}
We can then rewrite~\eqref{eqn:E_P-lower-semicontinuity1}, the inequality we seek to prove, as
\begin{align}
    \mathrm{e}(u_\infty)\leq \liminf_{k \to \infty} \mathrm{e}(u_k).
\end{align}

The issue we must address is that the Plateau solutions $u_k$ might not converge to $u_\infty$ as maps in any sense. Nevertheless, we will still be able to extract a weak subsequential limit $u_*$, whose Dirichlet energy, and therefore surface area, is equal to that of $u_\infty$.

\subsection{Equicontinuity at the Boundary}
We first prove the following simple lemma,
\begin{lemma}\label{lemma:diameter-convergence}
    Suppose $\left\{\mu_k:[a,b]\rightarrow\RR^3\right\}$ is a sequence of continuous curves that converges uniformly to $\mu_\infty$. Then the diameters of the images converge: $\mathrm{diam}(\mu_k([a,b]))\rightarrow\mathrm{diam}(\mu_\infty([a,b]))$.
\end{lemma}
\begin{proof}
    Let $\varepsilon>0$ and choose $N$ large such that $\norm{\mu_k-\mu_\infty}_{C^0}<\frac{\varepsilon}{2}$ for $k>N$. Let $x,y\in[a,b]$ such that $\abs{\mu_\infty(x)-\mu_\infty(y)}=\mathrm{diam}(\mu_\infty([a,b]))$. From uniform convergence we have, for any $s,t\in[a,b]$ and $k > N$,
    \begin{align*}
        \abs{\mu_k(s)-\mu_\infty(s)},~\abs{\mu_k(t)-\mu_\infty(t)}<\frac{\varepsilon}{2}.
    \end{align*}
    We have the following lower-bound for $\mathrm{diam}(\mu_k([a,b]))$,
    \begin{align*}
        \mathrm{diam}(\mu_\infty([a,b]))&=\abs{\mu_\infty(x)-\mu_\infty(y)}\\
        &\leq\abs{\mu_\infty(x)-\mu_k(x)}+\abs{\mu_k(x)-\mu_k(y)}+\abs{\mu_k(y)-\mu_\infty(y)}\\
        &\leq\mathrm{diam}(\mu_k([a,b]))+\varepsilon.
    \end{align*}
    Let $s,t\in[a,b]$, then we have
    \begin{align*}
        \abs{\mu_k(s)-\mu_k(t)}&\leq \abs{\mu_k(s)-\mu_\infty(s)}+\abs{\mu_\infty(s)-\mu_\infty(t)}+\abs{\mu_\infty(t)-\mu_k(t)}\\
        &\leq \mathrm{diam}(\mu_\infty([a,b]))+\varepsilon.
    \end{align*}
    Taking a $\max$ over $s,t$ in the above inequality yields
    \begin{align*}
        \mathrm{diam}(\mu_k([a,b]))\leq \mathrm{diam}(\mu_\infty([a,b]))+\varepsilon.
    \end{align*}
\end{proof}

We also state the Courant-Lebesgue lemma which is essential for proving equicontinuity of our maps at the boundary. First, for $p\in \bar D$ and $\rho>0$, we define the circles
\begin{align*}
    C_\rho (p) = \{q\in\bar D:~\abs{p-q}=\rho\}.
\end{align*}

\begin{lemma}[Courant-Lebesgue]\label{lemma:courant-lebesgue}
    Let $u:D\rightarrow\RR^3$ with $u\in C^0(\bar D)\cap W^{1,2}(D)$ where $\mathrm{e}(u)\leq K/2$ for some $K>0$. Then, for all $0 < \delta \ll 1$, there exist $\rho\in[\delta,\sqrt{\delta}]$ such that for all $p\in\bar D$,
    \begin{align*}
        \left[\mathrm{diam}(u(C_\rho(p)))\right]^2\leq\frac{8\pi^2 K}{-\log\delta}.
    \end{align*}
\end{lemma}
\begin{proof}
    See~\cite{colding2011course} Lemma 4.11.
\end{proof}

For the rest of this section, fix three distinct points $p_1,p_2,p_3\in\partial D$. For each $k$, we can find a unique conformal diffeomorphism $\phi_k:D\rightarrow D$ such that
\begin{align}\label{eqn:gauge-fixing}
    u_k\circ\phi_k(p_i)=\gamma_k(p_i)\text{ for }i=1,2,3.
\end{align}

\begin{proposition}\label{prop:equicontinuity}
    The maps $\tilde u_k:=u_k\circ\phi_k$ are equicontinuous on $\partial D$.
\end{proposition}
\begin{proof}
    Since $\{\tilde u_k\}$ is a minimizing sequence (and $\mathrm{e}[\cdot]$ is conformally invariant), we can assume that $\mathrm{e}(\tilde u_k)\leq K/2$ for some $K>0$. Let $\varepsilon>0$ and without loss of generality assume that
    \begin{align*}
        \varepsilon<\liminf_{k \to \infty}\left(\min_{i\neq j}\abs{\gamma_k(p_i)-\gamma_k(p_j)}\right),
    \end{align*}
    which exists since the $\gamma_k$ are all simple closed curves, along with $\gamma_\infty$. Since $\gamma_\infty$ is simple, closed and has finite arc-length, there is some $d_0>0$ such that for $p,q\in\partial D$ with $0<\abs{\gamma_\infty(p)-\gamma_\infty(q)}<d_0$, $\gamma_\infty(\partial D)\setminus\{p,q\}$ has exactly one component with diameter, $\mathrm{diam}_\infty\leq\frac{\varepsilon}{2}$~\cite[Proof of Lemma 4.14]{colding2011course}. Call the closure of this component $\Gamma$. Let $N>0$ such that $\norm{\gamma_k-\gamma_\infty}_{C^0}<\frac{\varepsilon}{4}$ for $k>N$. Then, from Lemma~\ref{lemma:diameter-convergence}, the sub-arc of $\gamma_k(\partial D)$ parametrized by the same part of $\partial D$ that parametrizes $\Gamma$ has diameter (denoted $\mathrm{diam}_k$) $\frac{\varepsilon}{2}$-close to the diameter of $\Gamma$, i.e. $\abs{\mathrm{diam}_k-\mathrm{diam}_\infty}<\frac\varepsilon 2$. Thus, $\mathrm{diam}_k<\varepsilon$. It follows that $\abs{\gamma_k(p)-\gamma_k(q)}<\frac{\varepsilon}{2}+d_0$. Henceforth, write $d:=\frac{\varepsilon}{2}+d_0$ and let $k>N$.
    
    Take $\delta<1$ small enough such that $\frac{8\pi^2 K}{-\log\delta}<d^2$ and not more than one of the $p_i$'s is in $B(p,\sqrt{\delta})$. Given any $p\in\partial D$, from Lemma~\ref{lemma:courant-lebesgue}, there exists $\rho\in[\delta,\sqrt{\delta}]$ such that
    \begin{align*}
        \textrm{diam}(\tilde u_k(C_\rho(p)))^2<\frac{8\pi^2 K}{-\log\delta}<d^2.
    \end{align*}
    In other words, $\textrm{diam}(\tilde u_k(C_\rho(p)))<d$. $C_\rho(p)$ divides $\partial D$ into two components, $A_1$ and $A_2$. Without loss of generality, let $A_1$ be the component that contains fewer than two of $p_1,p_2,p_3$. 
    
    Denote the images, $G_1= u(A_1)$ and $G_2= u(A_2)$. Due to monotonicity, $G_1$ will also contain fewer than two of $\gamma_k(p_1),\gamma_k(p_2),\gamma_k(p_3)$. Since $\textrm{diam}(\tilde u_k(C_\rho(p)))<d$, from Lemma~\ref{lemma:courant-lebesgue}, $\mathrm{diam}(G_i)<\varepsilon$ for at least one of $i=1,2$. However, from the chosen value of $\varepsilon$, this component cannot contain more than one of $\gamma_k(p_1),\gamma_k(p_2),$ or $\gamma_k(p_3)$. Thus, the component has to be $G_1$. 
    
    In conclusion, we have shown that there is $\rho>0$ such that $\abs{p-q}<\rho$ implies $\abs{\gamma_k(p)-\gamma_k(q)}<\varepsilon$.
\end{proof}

\subsection{Convergence in the Interior}
Since $\mathrm{e}(u_k)$ is uniformly bounded, we have from compactness that after passing to a subsequence, $u_k\rightharpoonup u_*$ weakly in $W^{1,2}(D,\RR^3)$ for some $u_*$, not necessarily equal to $u_\infty$, and 
\begin{align*}
    \mathrm{e}(u_*)\leq\liminf_{j \to \infty}\mathrm{e}(u_j).
\end{align*}
Furthermore, from Proposition~\ref{prop:equicontinuity}, $\{u_j\}$ is equicontinuous on $\partial D$ and by Arzela-Ascoli, a subsequence $\{u_j\}$ (not relabelled) converges uniformly on $\partial D$. Since each $u_j$ is harmonic, we have from the maximum principle that 
\begin{align*}
    \sup_D\abs{u_j-u_k}=\max_{\partial D}\abs{u_j-u_k}.
\end{align*}
Thus $\{u_j\}$ converges uniformly on $\bar{D}$ to $u_*$. Moreover, uniform convergence implies that $u_*:\partial D\rightarrow\gamma_\infty$ is monotone and onto, so therefore $u_*\in \mathcal{D}_{\gamma_\infty}$. We thus obtain the desired inequality,
\begin{align*}
    \mathrm{e}(u_\infty)=\inf_{u\in\mathcal{D}_{\gamma_\infty}}\mathrm{e}(u)\leq\mathrm{e}(u_*)\leq\liminf\limits_{j \to \infty}\mathrm{e}(u_j).
\end{align*}

With this, we have established the lower semicontinuity of $E$, and $\gamma_\infty\in\mathcal{A}$ is an energy minimizer. We can thus arrive at the main theorem for this section:
\begin{theorem}\label{thm:existence-of-minimiser}
    There exists $\gamma_\infty\in\mathcal{A}$ which minimizes $E$.
\end{theorem}

\section{Necessary Conditions of M\"obius-Plateau Critical Helix Pairs}\label{sec:stationary-helix-pairs}

In this section, we alter the problem slightly to consider the M\"obius energy of links, where $\gamma_1, \gamma_2: \R \to \R^3$ are disjoint locally rectifiable curves. The M\"obius energy of the link is defined as 
\begin{align}
\label{mobiuslink}
E_M(\gamma_1, \gamma_2) = \int_{\R \times \R} \frac{ |\dot{\gamma_1}(u)| | \dot{\gamma_2}(v) | dudv}{|\gamma_1(u) - \gamma_2(v)|^2}.
\end{align}
In the link case, we do not subtract an intrinsic distance term in the integrand, as it is not needed for the integral to converge. At a given point $\gamma_1(u)$, the $L^2$ gradient of the M\"obius energy is given by the vector-valued integral
\begin{align}
\label{mobiusvar}
    G_{\gamma_1, \gamma_2}(u) = 2 \int_\R \left[ \frac{2 P_{\dot{\gamma_1}(u)^\perp}(\gamma_2(v) - \gamma_1(u))}{|\gamma_2(v) - \gamma_1(u)|^2} - \textbf{N}_{\gamma_1(u)} \right] \frac{|\dot{\gamma_2}(v)|dv}{|\gamma_2(v) - \gamma_1(u)|^2}.
\end{align}
The derivation of the variational equation for the M\"obius energy of knots is given in \cite{freedman1994mobius}, and He \cite{he02} showed the equation also holds in the case of links. Here, $P_{\dot{\gamma_1}(u)^\perp}$ refers to the projection onto the plane normal to the tangent vector at $\gamma_1(u)$, and $\textbf{N}_{\gamma_1(u)}$ refers to the normal vector in the Frenet frame along $\gamma_1$. The gradient is defined similarly along $\gamma_2$ by switching $u$ and $v$.  

A rough heuristic for interpreting \eqref{mobiusvar} is that curves with large M\"obius gradient are those in which many rescaled pairwise difference vectors are close to Frenet binormals. The integral projects a rescaled difference vector onto the curve's normal plane, and subtracts the Frenet normal, leaving only the binormal. The rescaling obeys an inverse quartic law, so difference vectors between nearby points of the curve which approximate the binormal well will contribute greatly to the gradient integral. The binormal of a helix points roughly in the direction of its axis, hinting that helices provide a fruitful set of examples in studying the dynamics of the M\"obius energy. Using explicit helix parametrizations turns \eqref{mobiusvar} into an improper trigonometric integral, which we analyze with power series and elementary estimates. In a forthcoming paper, the first author analyzes the M\"obius variations for helix pairs using complex variable techniques distinct from the methods in this section.

Whilst both energies are infinite, the variational equation describing critical parametrizations remains well-defined. Let $\boldsymbol\nu$ be the oriented unit surface normal to $\Sigma$. Now, the direction of the $L^2$ gradient of the Plateau energy at $\gamma_1(u)$ (resp. $\gamma_2(v)$) is the unit conormal $\textbf{n} = -\textbf{T} \times \boldsymbol\nu$ (resp. $\textbf{T} \times \boldsymbol\nu$). This unit vector is oriented so that it is pointing away from the surface at each of the boundary components, which in turn defines the orientation of $\boldsymbol\nu$. The derivation of the gradient for the Plateau energy is a standard calculation which can be found in Prop. 1 of \cite{gruber21}. Thus, the variational equation for the M\"obius-Plateau energy is
\begin{align}
\label{mobiusplateauvar}
    \alpha G_{\gamma_1,\gamma_2} = - \beta \bf{n}.
\end{align}

A general helicoid $\Sigma = \Sigma(\omega)$ has parametrization 
\begin{align}
\label{helicoid}
    \begin{bmatrix}
    x(s,t) \\
    y(s,t) \\
    z(s,t)
    \end{bmatrix} &=
    \begin{bmatrix}
    s \cos(\omega t) \\
    s \sin(\omega t) \\
    t
    \end{bmatrix}.
\end{align}

\begin{definition}
A \textbf{helicoidal strip} is the parametrized surface in \eqref{helicoid} with $A \leq s \leq B, t \in \R$, for prescribed $A < B$. A \textbf{helicoidal screw} corresponds to the case where $A$ and $B$ have opposite signs, and a \textbf{helicoidal ribbon} corresponds to the case where $A$ and $B$ have the same signs.
\end{definition}
The two boundary curves $\gamma_1$ and $\gamma_2$ of a helicoidal strip are found by setting $s = A$ and $s = B$. We will parametrize the curves with $u$ and $v$ respectively. 

\begin{remark}
    In the prior sections, we were careful in considering the possibility of sudden transitions in the Plateau solutions induced by small variations of the boundary curve. As the frequency parameter of a helix changes, there can be a so-called ``Courant headphone transition" where the area-minimizing surface jumps from a helicoid to a surface between the vertical layers resembling a telephone cord. In Courant's time, headphone cords probably looked like this, hence the name. The transition parameters are discussed in \cite{courant40} and \cite{machon16}. As this section is concerned with necessary conditions for critical helix pairs, we do not take these transitions into account. As we will see later, numerical searches for critical helix pairs bounding a ribbon show their widths are orders of magnitude smaller than the frequency, precluding the possibility of a headphone transition. In other words, the Plateau solution for the critical helix pair is their common helicoid.
\end{remark}

The space of helicoidal screws and ribbons are given by the three parameters $A,B,$ and $\omega$. A helicoidal strip is a screw if and only if it contains the helical axis, and otherwise it is a ribbon. For given $\alpha$ and $\beta$, we seek to find the values of $A,B,$ and $\omega$ which satisfy the variational equation of the M\"obius-Plateau energy. The variational equations for stationary helix pairs differs in the case of screws and ribbons, because of the differing orientations of $\mathbf{N}$, which we will see results in two trigonometric integrals that nominally look similar, but greatly vary in how the signs of their evaluations depend on their parameters.

For fixed $\omega$, the boundary helices are invariant under ``screw" transformations which rotate the $xy$-plane by angle $\omega t$ whilst translating in the $z$ direction by $t$. As the screw transformations are a subgroup of Euclidean isometries (and also of the larger M\"obius group of $\R^3$), the gradient flow of the M\"obius energy starting at a helix preserves symmetries under the screw transformations (cf. p.42 of \cite{freedman1994mobius}). Therefore, $G_{\gamma_1,\gamma_2}$ is tangent to $\Sigma(\omega)$, and to calculate the entire M\"obius gradient vector fields along $\gamma_1$ and $\gamma_2$, it suffices to compute them at two particular points along the boundary curves and then apply screw transformations. The same symmetry rule also holds for computing the Plateau gradient vector fields along the helices, as the Plateau energy is invariant under Euclidean isometries.

Observe that 
\begin{align*}
    \dot{\gamma_1}(u) &= 
    \begin{bmatrix}
    -\omega A \sin(\omega u) \\
    \omega A \cos(\omega u) \\
    1
    \end{bmatrix},
\end{align*} and that $|\dot{\gamma_1}(u)| = \sqrt{\omega^2 A^2 + 1}$ for all $u$. We also have $|\dot{\gamma_2}(u)| = \sqrt{\omega^2 B^2 + 1}$ for all $v$. Furthermore,
\begin{align*}
    \Ddot{\gamma_1}(u) = \begin{bmatrix}
    -\omega^2 A \cos(\omega u) \\
    -\omega^2 A \sin(\omega u) \\
    0
    \end{bmatrix}.
\end{align*}

Next, for given surface parameters $s$ and $t$, the tangent plane to $\Sigma$ is given by 
\begin{align*}
    T\Sigma &= \text{span} \Biggl\{ \partial_s \begin{bmatrix}
    x \\
    y \\
    z
    \end{bmatrix}, \partial_t \begin{bmatrix}
    x \\ 
    y \\ 
    z 
    \end{bmatrix} \Biggr\} \\
    &= \text{span} \Biggl\{ \begin{bmatrix}
    \cos(\omega t) \\
    \sin(\omega t) \\ 
    0
    \end{bmatrix}, \begin{bmatrix}
    -\omega s \sin(\omega t) \\
    \omega s \cos(\omega t) \\
    1
    \end{bmatrix} \Biggr\}.
\end{align*}

So at $u = 0$, we have the following:
\begin{center}
    $\gamma_1(0) = \begin{bmatrix}
    A\\
    0\\
    0
    \end{bmatrix}, \dot{\gamma_1}(0) = \begin{bmatrix}
    0\\
    \omega A \\
    1
    \end{bmatrix}, \Ddot{\gamma_1}(0) = \begin{bmatrix}
    -\omega^2 A\\
    0\\
    0
    \end{bmatrix}$.
\end{center}

From this, we can see that $\textbf{N}_{\gamma_1(0)} = (\pm 1, 0, 0)$, depending on the sign of $A$. We can also see that $T_{\gamma_1(0)} \Sigma = \text{span} \{ (1, 0, 0), (0, \omega A, 1) \}$, and thus $\boldsymbol\nu = \frac{1}{\sqrt{\omega^2A^2 + 1}}(0, \omega A , 1)$. Hence, $\textbf{n} = (-1,0,0)$.

\begin{figure}
    \centering
    \begin{subfigure}{0.5\textwidth}
    \centering
    \includegraphics[width=4in]{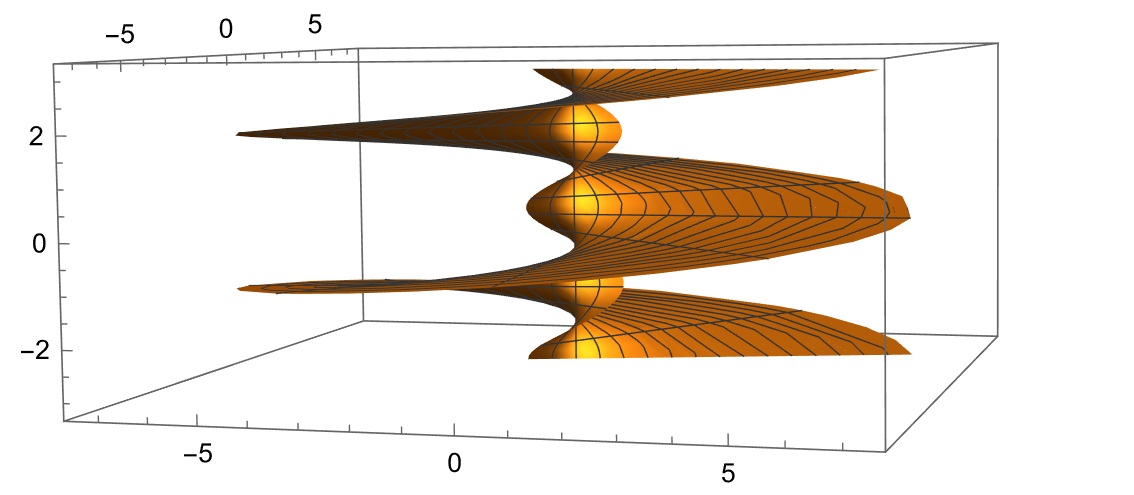}
    \caption{A stationary asymmetric helicoidal screw.}
    \end{subfigure}%
    \begin{subfigure}{0.5\textwidth}
    \centering
    \includegraphics[width=1in]{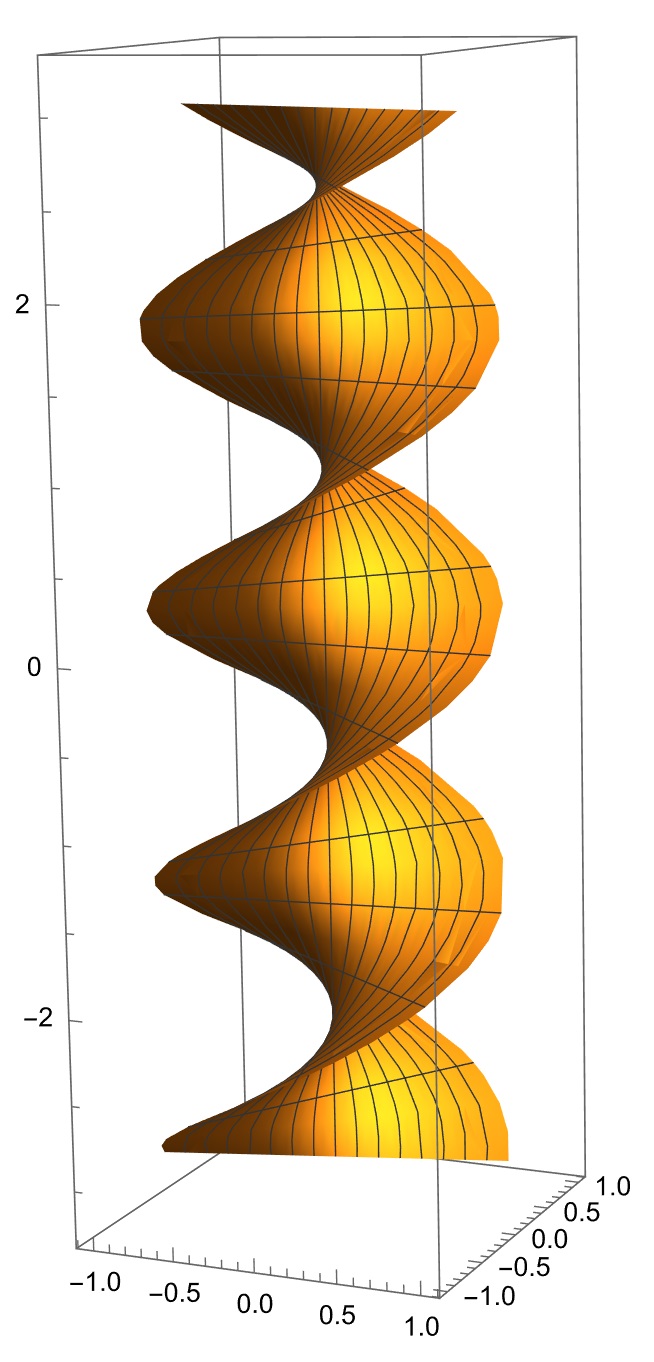}
    \caption{A stationary symmetric helicoidal screw.}
    \end{subfigure}
    \caption{Two stationary helicoidal screws with identical parameters $A = -1, \omega = 2$, except $B \approx 6.15375$ in the figure on the left, whilst $B = 1$ in the figure on the right. In these examples, we set $\alpha = 2, \beta = 1$ }
    \label{stablescrews}
\end{figure}

Next, observe
\begin{align*}
    P_{\dot{\gamma}_1(0)^\perp} \big(\gamma_2(v) - \gamma_1(0)\big) &= \big(\gamma_2(v) - \gamma_1(0)\big) - \langle \gamma_2(v) - \gamma_1(0), \dot{\gamma_1}(0) \rangle \frac{\dot{\gamma_1}(0)}{|\dot{\gamma_1}(0)|^2} \\
    &=  \begin{bmatrix}
    B \cos(\omega v) - A \\
    B \sin(\omega v) \\
    v
    \end{bmatrix} -  \frac{ \omega AB\sin(\omega v) + v}{\omega^2A^2 + 1} \begin{bmatrix}
    0\\
    \omega A \\
    1
    \end{bmatrix}\\
    &= \begin{bmatrix}
    B\cos(\omega v) - A\\
    \frac{B \sin(\omega v) - A \omega v}{\omega^2 A^2 + 1} \\
    v - \frac{v + \omega AB \sin(\omega v)}{\omega^2A^2 + 1}
    \end{bmatrix}.
\end{align*}
Finally, it is straightforward to see that $|\gamma_2(v) - \gamma_1(0)|^2 = A^2 - 2AB\cos(\omega v) + B^2 + v^2$. For now, we will make the assumption that $A < 0 < B$, which means $\textbf{N} = (1,0,0)$. Putting all of this together, the variational equation becomes the system:
\begin{align}
    2 \alpha \int_{-\infty}^\infty \left[ \frac{2 \left( B \cos(\omega v) - A \right)}{A^2 - 2AB\cos(\omega v) + B^2 + v^2} - 1 \right] \frac{\sqrt{\omega^2 B^2 + 1}}{A^2 - 2AB\cos(\omega v) + B^2 + v^2} dv &= \beta \label{gamma1var}\\
    \int_{-\infty}^\infty \frac{B \sin(\omega v) - \frac{v + \omega^2 A^2 B \sin(\omega v)}{\omega^2 A^2 + 1}}{(A^2 - 2AB\cos(\omega v) + B^2 + v^2)^2}dv &= 0 \nonumber\\
    \int_{-\infty}^\infty \frac{ v - \frac{\omega A B \sin(\omega v) + v}{\omega^2 A^2 + 1}}{\left( A^2 - 2AB\cos(\omega v) + B^2 + v^2 \right)^2} dv &= 0.\nonumber
\end{align}
The last two integrands are odd functions in $v$, which means the integrals are always going to be zero, so only equation \eqref{gamma1var} relevant to us. Through identical computations by taking $v = 0$, noting that in this case we have $\textbf{N}_{\gamma_2(0)} = (-1,0,0)$ and $\textbf{n} = (1, 0, 0)$, we see the first component of the variational equation at $\gamma_2(0)$ is
\begin{align}
     2 \alpha \int_{-\infty}^\infty \left[ \frac{2 \left( A \cos(\omega u) - B \right)}{A^2 - 2AB\cos(\omega u) + B^2 + u^2} + 1 \right] \frac{\sqrt{\omega^2 A^2 + 1}}{A^2 - 2AB\cos(\omega u) + B^2 + u^2} du &= -\beta. \label{gamma2var}
\end{align}
Relabelling the variable of integration and adding \eqref{gamma1var} and \eqref{gamma2var} yields
\begin{equation}
\label{screwvar}
    4 \alpha \int_{-\infty}^\infty \frac{(A + B)(\cos(\omega v) - 1)}{\left(A^2 - 2AB \cos(\omega v) + B^2 + v^2 \right)^2}dv = \beta \left( \frac{1}{\sqrt{\omega^2 B^2 + 1}} - \frac{1}{\sqrt{\omega^2 A^2 + 1}}\right).
\end{equation}
This single equation leaves us with a two parameter family of stationary helicoids controlled by the constants $\alpha$ and $\beta$. From direct inspection, it is clear that $A = -B$ is a solution, which we call the symmetric stationary helicoidal screw. It is straightforward to numerically search for other solutions of \eqref{screwvar}. For instance, setting $\alpha = 2, \beta = 1, A = -1, \omega = 2$, and solving for $B$ gives solutions $B = 1$ and $B \approx 6.15375$. The two helicoidal screws are pictured in \ref{stablescrews}. In the numerical searches we ran, we generally found that under most choices of parameters, the symmetric solution had a larger basin of attraction than the asymmetric solution. An analysis of the attraction properties of the solutions is yet to be done.

The full description of this family, especially with respect to its singularities, is not fully known. When we tried to find solutions to the variational equation for helicoidal ribbons, derived below, we had difficulties finding solutions, and we suspected at first that no solutions existed. During our attempt to prove this lack of solutions, we had difficulty in bounding the negative contribution to the variational integral, and we later found that solutions do indeed exist, but only under comparatively strict conditions. Furthermore, the solutions we did manage to find have parameters which differ in orders of magnitude, which contributed to the difficulty in the numerical search.

\begin{theorem}
A stationary helicoidal ribbon must have high frequency (equivalently, low pitch), very small width in comparison to the frequency, and remain close to the axis. That is, the following conditions on the parameters must hold:
\begin{enumerate}
    \item $\omega > 1$
    \item $|B - A|^2 < C(\alpha,\beta) \omega$
    \item $|A + B| < 2$.
\end{enumerate}
\end{theorem}
\begin{remark}
    We only prove the case when $\alpha = \beta = 1$. In the general case, the bounds \eqref{ba4w} and \eqref{cosneg} are scaled by powers of $\alpha$ and $\beta$, though \eqref{tail} is unchanged. One may verify these bounds still sum to a positive number as required in the last step of the proof, but this only further complicates the algebra without fundamentally altering our reasoning. Our numerical search suggests the second condition could be considerably strengthened to the form $|B-A|^{k_1} < C \frac{\beta^{k_2}}{\alpha^{k_3}} \omega$, with the optimal exponents $k_i$ yet to be determined. Modifying the bound in \eqref{cosneg} is more difficult since we could get a family of rational functions in $\omega$ controlled by $\alpha$ and $\beta$, and so a numerical approximation does not rigorously justify a lower bound. When applying the M\"obius-Plateau energy to mainstream Euler-Plateau problems, it is likely necessary to vary the coefficients $\alpha$ and $\beta$, so understanding their precise scaling effects is imperative.
\end{remark}
\begin{proof}
Without loss of generality assume $0 < A < B$. We may also assume $\omega > 0$ because the cosine function is even. This is tantamount to fixing the chirality of the helicoids under consideration to being right-handed. Observe $\textbf{N}_{\gamma_1(0)} = (-1,0,0)$, so the $-1$ term in \eqref{gamma1var} corresponding to $-\textbf{N}$ becomes a $+1$ whilst \eqref{gamma2var} remains unchanged. Adding these two nontrivial variational equations at $\gamma_1(0)$ and $\gamma_2(0)$ yields
\begin{align}
\label{ribbonvar}
    4 \int_{-\infty}^\infty \frac{A^2 - A + B^2 - B + (A + B - 2AB)\cos(\omega v) + v^2}{\left(A^2 - 2AB \cos(\omega v) + B^2 + v^2 \right)^2}dv &= \left( \frac{1}{\sqrt{\omega^2 B^2 + 1}} - \frac{1}{\sqrt{\omega^2 A^2 + 1}}\right).
\end{align}
We will show that except in a very limited circumstance, the integrand on the lefthand side of \eqref{ribbonvar} is always positive. When the integrand can take negative values, we will then show that unless all of the conditions listed are met, then the integral is positive. The righthand side of \eqref{ribbonvar} is negative, as $B > A$, so we will get our desired result.

The denominator is obviously positive. As for the numerator, we break down into cases. Suppose $A + B  -2AB \leq 0$. Then
\begin{align*}
    A^2 - A + B^2 - B + (A + B - 2AB)\cos(\omega v) + v^2 &\geq A^2 - A + B^2 + A + B - 2AB + v^2\\
    &= (B - A)^2 + v^2 > 0. 
\end{align*}
For the remainder of the proof, assume $A + B - 2AB > 0$ and we will break down into further cases. Assume $A + B \geq 2$. Then
\begin{align*}
  A^2 - A + B^2 - B + (A + B - 2AB)\cos(\omega v) + v^2 &\geq A^2 - A + B^2 - B -(A + B - 2AB) + v^2\\
    &= (A + B)(A + B - 2) + v^2 > 0. 
\end{align*}
Next, assume $ A + B < 2$ with $\omega \leq 1$. From the Taylor formula, $\cos(\omega v) \geq 1 - \frac{\omega^2v^2}{2} \geq 1 - \frac{v^2}{2}$. Now see that
\begin{align*}
    A^2 - A + B^2 - B + (A + B - 2AB)\cos(\omega v) + v^2 &\geq A^2 - A + B^2 - B + (A + B - 2AB)\left(1 - \frac{v^2}{2}\right) + v^2 \\
    &= (B - A)^2 + \frac{v^2}{2}(2AB - A - B + 2) \\
    &> \frac{v^2}{2}(2AB - A - B + 2)\\
    &\geq \frac{v^2}{2}(2AB) \geq 0.
\end{align*}
So in all of the cases we have considered for $A, B$, and $\omega$ so far, the integrand in \eqref{ribbonvar} is positive and so the variational equation will have no solutions.

\begin{figure}
\begin{center}
\includegraphics[width=4in]{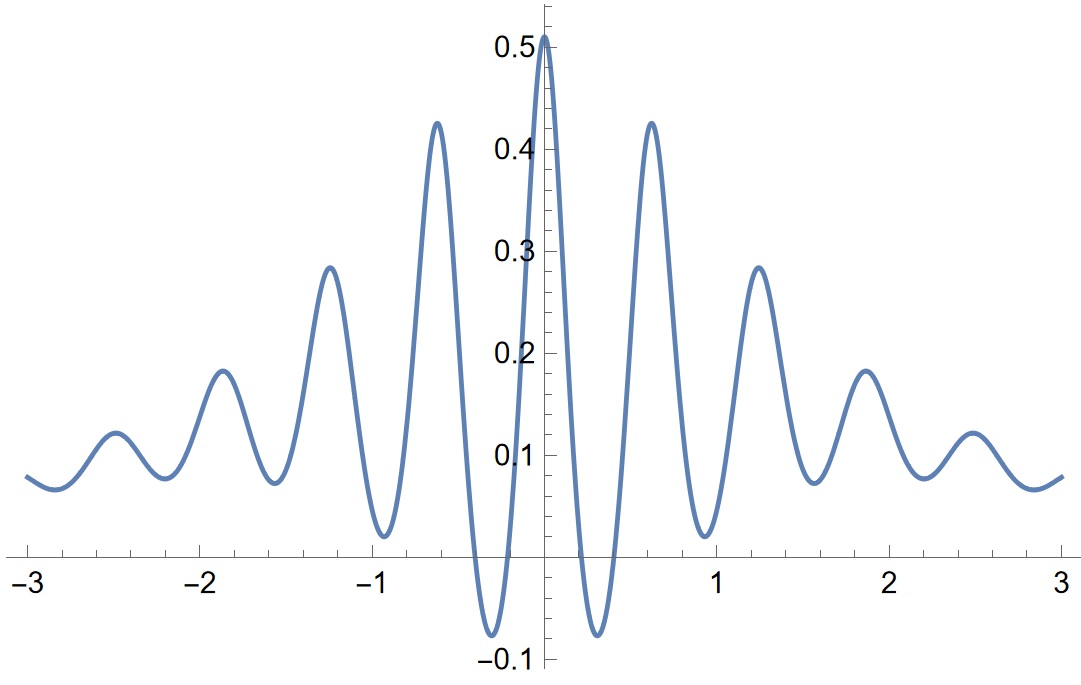}
\caption{A graph of the integrand in \eqref{ribbonvar} with $A = 0.1, B = 1,$ and $\omega = 10$, as a function of $v$.}
\label{ribleftgraph}
\end{center}
\end{figure}

Now, assume $A + B < 2$ with $\omega > 1$. It is only in this scenario that the numerator could take negative values and possibly lead to the integral in \eqref{ribbonvar} being negative. For instance, when $A = 0.1, B = 1, \omega = 10$, the integrand takes negative values as shown in Figure~\ref{ribleftgraph}.

As the integrand is an even function in $v$, it suffices to prove 
\begin{align*}
    \int_0^\infty \frac{A^2 - A + B^2 - B + (A + B - 2AB)\cos(\omega v) + v^2}{(A^2 - 2AB\cos(\omega v) + B^2 + v^2)^2}dv > 0,
\end{align*}
when the second condition in our theorem is negated. Our method will be to overestimate the magnitude of the contribution of the negative part of the integrand whilst underestimating the contribution of the positive part of the integrand, and show that the positive contribution still outweighs the negative contribution. To estimate the negative contribution, we integrate the quadratic Maclaurin series over $[0,\frac{1}{4}]$, and then compute a lower Riemann sum over  $\{\cos(\omega v) < 0\} \cap [0,1]$. When $(B-A)^2 > c\omega$, where $c$ is some constant independent of the parameters, we get \eqref{ba4w} and \eqref{cosneg}. The computation is detailed in an appendix.

As we saw, for $v \geq 1$, the integrand is strictly positive. So consider the contribution of the tail
\begin{align*}
    \int_1^\infty \frac{A^2 - A + B^2 - B + (A + B - 2AB)\cos(\omega v) + v^2}{(A^2 - 2AB\cos(\omega v) + B^2 + v^2)^2}dv &> 0.
\end{align*}
We see that $A^2 - 2AB\cos(\omega v) + B^2 + v^2 \geq (A + B)^2 + v^2$. Combining this with \eqref{supab4}, we can conclude that the denominator is bounded above by $(4 + v^2)^2$. Likewise, we again observe by \eqref{numbound} that the numerator is bounded below by $(A + B)(A + B - 2) + v^2$. Therefore, the numerator is bounded below by $v^2 - 1$. Putting all this together, we have
\begin{align}
    \int_1^\infty \frac{A^2 - A + B^2 - B + (A + B - 2AB)\cos(\omega v) + v^2}{(A^2 - 2AB\cos(\omega v) + B^2 + v^2)^2}dv &\geq \int_1^\infty \frac{v^2 - 1}{(4 + v^2)^2}dv \nonumber \\
    &= \frac{1}{32} \left(4 + 3\pi - 6\text{arccot}(2) \right) \nonumber\\
    &\geq \frac{3}{10} \label{tail}.
\end{align}
By adding our bounds from \eqref{ba4w},\eqref{cosneg}, and \eqref{tail}, we can conclude the integral is positive.
\end{proof}

\begin{figure}
    \centering
    \includegraphics[width=3in]{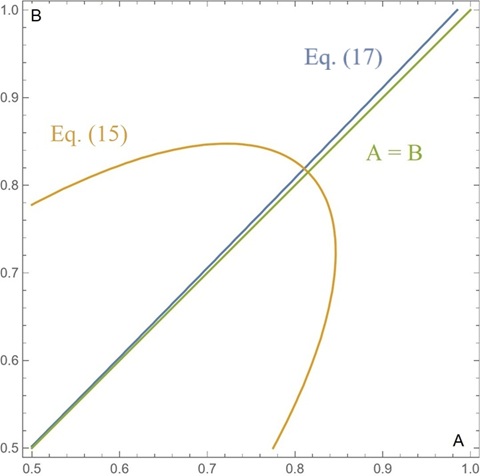}
    \caption{Contours in the $AB$-plane for \eqref{ribbonvar} and \eqref{ribbondiff}, with $\alpha = \beta = 1, \omega = 10$, and approximate region of integration $[-10000,10000]$. The intersection of the contour curves is $(A,B) \approx (0.810515,0.818991)$ which is in the valid parameter space $\{(A,B): 0 < A < B, A + B < 2 \}$.}
    \label{ribboncontours}
\end{figure}

Negating the second condition was crucial in our proof, as we otherwise would not be able to obtain the bound in \eqref{ba4w}, as the coefficient $-\frac{1}{(B-A)^4}$ tends to $-\infty$ as $(B-A) \to 0^+$, which makes the bound useless in limiting the contribution of the negative parts of the integrand. We also used this assumption in \eqref{cosnegden} in order to force the Riemann sum approximation to be of order $O(\omega^{-4})$. Were we to na\"ively the use the inequality $A^2 + B^2 + v^2 \geq 0 + \left(\frac{\pi}{2\omega}\right)^2$, our Riemann sum approximation would be of order $O(\omega^4)$ with a leading negative coefficient, which tends to $-\infty$ as $\omega \to \infty$, again preventing a bound for the negative contribution.

With the restrictions on the parameters in mind, we found approximate numerical solutions to the variational equations for stationary helicoidal ribbons. Our analysis so far has only dealt with the sum of the variational equations for $\gamma_1$ and $\gamma_2$, so to solve the full variational system, we included the difference of the two variational equations, which is
\begin{align}
    4 \alpha \int_{-\infty}^\infty \frac{(B - A) \cos^2\left( \frac{\omega v}{2} \right)}{\left(A^2 - 2AB\cos(\omega v) + B^2 + v^2 \right)^2}dv &= \beta \left( \frac{1}{\sqrt{\omega^2 A^2 + 1}} + \frac{1}{\sqrt{\omega^2 B^2 + 1}} \right).
    \label{ribbondiff}
\end{align}
Both sides of this equation are clearly positive for $B > A$, and hence there are no further obvious restrictions on the parameters for $A$ and $B$, which is why we opted to study \eqref{ribbonvar} in order to characterize the stationary ribbons.

In our numerical search, we evaluated the integrals corresponding to \eqref{ribbonvar} and \eqref{ribbondiff} on the finite interval $[-10000,10000]$ via a Gaussian quadrature with precision set to $10^{-6}$. With $\alpha = \beta = 1$ and $\omega = 10$, we obtained the contour plot depicted in \ref{ribboncontours}. The intersection of the contour lines for \eqref{ribbonvar} and \eqref{ribbondiff} was found to be $(A,B) \approx (0.810515,0.818991)$, which corresponds to a width of $B - A \approx 0.008476$. The effects of this intersection point under bifurcations of $\alpha, \beta$, and $\omega$ is yet to be investigated.

\section{Acknowledgements}
The authors would like to thank Xin Zhou for his numerous helpful discussions about this work, in particular making the initial suggestion to investigate the M\"obius-Plateau energy. The authors would also like to thank Steven Strogatz and Tobias Colding for their continued support.
\\
\\
\noindent This work is partially funded by an NSF RTG grant titled Dynamics, Probability and PDEs in Pure and Applied Mathematics, DMS-1645643. Max Lipton is supported by an NSF Postdoctoral Fellowship, DMS-2303384.
\appendix
\section{Estimating the negative contribution}
Suppose $\cos(\omega v) \geq 0$, and since by assumption $A+B-2AB>0$, we have
\begin{align*}
    A^2 - A + B^2 - B + (A + B - 2AB)\cos(\omega v) + v^2 &\geq A^2 - A + B^2 - B + v^2 \\
    &\geq -\frac{1}{2} + v^2.
\end{align*}
The last inequality follows from the fact that
\begin{align}
\label{infminusonehalf}
    \inf_{\substack{A + B < 2 \\ 0 < A < B}} A^2 - A + B^2 - B = -\frac{1}{2}.
\end{align}
Therefore, a necessary condition for the integrand to be negative whilst $\cos(\omega v) > 0$ is for $v \leq \frac{1}{4}$.

With a direct calculation or with a computer algebra system, one obtains the Macularin series of the integrand in \eqref{ribbonvar} expanded at $v = 0$:
\begin{align}
    \frac{A^2 - A + B^2 - B + (A + B - 2AB)\cos(\omega v) + v^2}{(A^2 - 2AB\cos(\omega v) + B^2 + v^2)^2} &= \frac{1}{(B-A)^2} - \frac{2 + (A + B + 2AB)\omega^2}{2(B-A)^4}v^2 + O(v^4).
    \label{maclaurin}
\end{align}
 Writing $f$ in place of the lefthand side of \eqref{maclaurin} and differentiating with respect to $v$, we can apply Fa\'a di Bruno's formula to see that 
\begin{align*}
    | f^{(k)}(0) | &\leq \left( \frac{c \omega}{(B-A)^2} \right)^k k!,
\end{align*}
where $c$ is a constant independent of $A,B,$ and $\omega$. Estimating $c$ independently requires the prior assumption that $|A+B| < 2$. So when $f$ is extended meromorphically to $\C$, the Taylor series converges for all $\{z: |z| \leq \delta\}$, where $\delta < \frac{(B-A)^2}{c\omega}$. By the maximum principle and the fact $f$ and its Taylor polynomials are even, we can also conclude that the differences between $f$ and its Taylor polynomials are strictly positive or negative for all $|v| \leq \delta$. Our argument requires $\delta = \frac{1}{4}$ to be a valid choice, which holds when $c \omega < (B-A)^2$, after a rescaling of $c$. Now assume this.

The coefficient of the $v^4$ term, which we omit, is rational function in $A,B,$ and $\omega$ with positive coefficients. Therefore,
\begin{align*}
    \frac{A^2 - A + B^2 - B + (A + B - 2AB)\cos(\omega v) + v^2}{(A^2 - 2AB\cos(\omega v) + B^2 + v^2)^2} &\geq \frac{1}{(B-A)^2} - \frac{2 + (A + B + 2AB)\omega^2}{2(B-A)^4}v^2.
\end{align*}
We will now integrate the quadratic Macluarin polynomial over $[0,\frac{1}{4}]$ to get a lower bound of the integral over all the points where the integrand could take negative values despite $\cos(\omega v) > 0$. Observe
\begin{align}
    &\int_0^{\frac{1}{4}} \frac{A^2 - A + B^2 - B + (A + B - 2AB)\cos(\omega v) + v^2}{(A^2 - 2AB\cos(\omega v) + B^2 + v^2)^2} dv \nonumber \\
    &\geq \int_0^{\frac{1}{4}} \left[\frac{1}{(B-A)^2} - \frac{2 + (A + B + 2AB)\omega^2}{2(B-A)^4}v^2 \right] dv \nonumber \\
    &= \frac{1}{4(B-A)^2} - \frac{2 + (A + B + 2AB)\omega^2}{384(B-A)^4} \nonumber \\
    &\geq  - \frac{2 + (A + B + 2AB)\omega^2}{384(B-A)^4} \nonumber \\
    &\geq -\frac{1 + \omega^2}{192(B-A)^4} \nonumber \\
    &> -\frac{1}{96}. \label{ba4w}
\end{align}
In the third to last inequality, we use the fact that
\begin{align}
    \sup_{\substack{A + B < 2 \\ 0 < A < B}} A + B - 2AB &= 2. \label{supab4}
\end{align}
Furthermore, in the last two inequalities, we use the assumption that $(B - A) \geq \omega > 1$.

Notice that outside of the domain of integration in \eqref{ba4w}, we will have $v \geq \frac{1}{4}$, and so the integrand is strictly positive when $\cos(\omega v) > 0$. In the remainder of the domain of integration, the integrand can only be negative provided that $\cos(\omega v) < 0$. In this instance, the numerator is bounded below by 
\begin{align}
    A^2 - A + B^2 - B + (A + B - 2AB)\cos(\omega v) + v^2 &\geq A^2 - A + B^2 - B - (A + B - 2AB) + v^2 \nonumber \\
    &= (A + B)(A + B - 2) + v^2 \label{numbound}.
\end{align}
Next, it is a straightforward optimization exercise to show that
\begin{align*}
    \inf_{\substack{A + B < 2\\ 0 < A < B}} (A + B)(A + B - 2) &= -1.
\end{align*}
Even though this infimum cannot be attained in the triangular domain of $A$ and $B$, if we fix any choice of $A$ and $B$, the quantity $A^2 - A + B^2 - B + (A + B - 2AB)\cos(\omega v)$ attains its largest negative magnitude when $v = \frac{n \pi}{ \omega}$, with $n$ an odd integer, and this quantity is bounded below by $-1$ over all valid $A$ and $B$. However, we have to add the positive term $v^2$ to get the total numerator, and so it is necessary for $v^2 < 1$ in order for the total numerator, and thus the integrand, to be negative.

The domains for $v$ such that $\cos(\omega v) \leq 0$ are the intervals $\left[\frac{n \pi}{\omega} - \frac{\pi}{2\omega}, \frac{n\pi}{\omega} + \frac{\pi}{2\omega}\right]$, with $n$ an odd integer. Each interval has length $\frac{\pi}{\omega}$. However, we also require $\left( \frac{n\pi}{\omega} - \frac{\pi}{2\omega} \right)^2 < 1$, . There are only finitely many of these intervals, which is the number of positive odd integers $n$ satisfying $n < \frac{1}{2} + \frac{\omega}{\pi}$. Hence, there are no more than $\frac{1}{4} + \frac{\omega}{2 \pi}$ such intervals. On the $n^{\text{th}}$ interval, the positive $v^2$ term is bounded below by $\left(\frac{(2n+1)\pi}{2\omega} - \frac{\pi}{2\omega}\right)^2$, and $\left(\frac{\pi}{2\omega}\right)^2$ independent of $n$. Hence to bound the denominator's magnitude from below to overestimate the negative contribution, observe 
\begin{align}
    A^2 - 2AB \cos(\omega v) + B^2 + v^2 &\geq A^2 + B^2 + v^2 \nonumber \\
    &= (B - A)^2 + 2AB + v^2 \nonumber \\
    &\geq \omega^2 + \left(\frac{\pi}{2\omega}\right)^2. \label{cosnegden}
\end{align}
Thus the magnitude of the denominator is bounded below by $\omega^4 + \frac{\pi^2}{2} + \frac{\pi^4}{16 \omega^4}$. Therefore, the negative contribution to the integral when $\cos(\omega v) < 0$ is bounded below by the Riemann sum
\begin{align}
\left( \frac{\pi}{\omega} \right) \sum\limits_{n = 0}^{\lceil \frac{1}{4} + \frac{\omega}{2\pi} \rceil} \left[ \frac{-1 + \left[ \frac{(2n + 1) \pi}{\omega} - \frac{\pi}{2\omega} \right]^2}{\omega^4 + \frac{\pi^2}{2} + \frac{\pi^4}{16 \omega^4}} \right] &= \left( \frac{\pi}{\omega^5 + \frac{\pi^2 \omega}{2} + \frac{\pi^4}{16 \omega^3}} \right) \sum\limits_{n = 0}^{\lceil \frac{1}{4} + \frac{\omega}{2\pi} \rceil} \left[ -1 + \frac{\pi^2}{4\omega^2} + \frac{2\pi^2}{\omega^2}n + \frac{4\pi^2}{\omega^2}n^2 \right]\nonumber\\
&\geq \left( \frac{\pi}{\omega^5 + \frac{\pi^2 \omega}{2} + \frac{\pi^4}{16 \omega^3}} \right) \left( \frac{1}{4} + \frac{\omega}{2 \pi} \right) \left( -1 + \frac{\pi^2}{4\omega^2} \right) \nonumber\\
&\hspace{0.5in}+ \left( \frac{2\pi^3}{\omega^7 + \frac{\pi^2 \omega^3}{2} + \frac{\pi^4}{16 \omega}} \right)\sum\limits_{n = 0}^{\lceil \frac{1}{4} + \frac{\omega}{2\pi} \rceil} \left(n + 2n^2 \right) \nonumber \\
&\geq  \frac{\omega (\pi - 2 \omega)(\pi + 2\omega)^2}{\left( \pi^2 + 4\omega^4 \right)^2}\nonumber\\
&\hspace{0.5in}+ \left( \frac{2\pi^3}{\omega^7 + \frac{\pi^2 \omega^3}{2} + \frac{\pi^4}{16 \omega}} \right)\Biggl[\frac{1}{2}\left(\frac{1}{4} + \frac{\omega}{2\pi}\right)\left(\frac{5}{4} + \frac{\omega}{2\pi} \right)  \nonumber \\
&\hspace{0.5in}+ \frac{1}{3}\left(\frac{1}{4} + \frac{\omega}{2\pi}\right)\left(\frac{5}{4} + \frac{\omega}{2\pi} \right)\left(\frac{3}{2} + \frac{\omega}{\pi}\right)\Biggr] \nonumber\\
&= \frac{\omega (\pi - 2 \omega)(\pi + 2\omega)^2}{\left( \pi^2 + 4\omega^4 \right)^2} + \frac{2\omega(3\pi + \omega)(\pi + 2 \omega)(5\pi + 2\omega)}{3\left(\pi^2 + 4\omega^4\right)^2} \nonumber\\
&= \frac{\omega(\pi + 2\omega)(33\pi^2 + 22\pi \omega - 8\omega^2)}{3\left(\pi^2 + 4 \omega^4\right)^2}.
\label{cosneg}
\end{align}
This lower bound is a rational function in $\omega$ of degree $-4$ dominated by the term $-\frac{1}{3\omega^4}$ in its Laurent expansion about $\frac{1}{\omega}$. It is straightforward to verify that \eqref{cosneg} has a global minimum over all $\omega > 1$ of approximately $-1.77 \times 10^{-6}$ at $\omega \approx 14.8$.

\bibliographystyle{amsalpha}
\bibliography{ref}

\end{document}